\documentclass[cmp,referee,envcountsame]{svjour}
\usepackage{amsmath}
\usepackage{amssymb}

\newcommand\bur{\operatorname{span}}

\newcommand\la{\langle}
\renewcommand\L{\mathcal L}
\newcommand\N{\mathbb N}
\newcommand\ot{\otimes}
\renewcommand\P{\mathcal P}
\newcommand\ra{\rangle}
\newcommand\rng{\operatorname{rng}}
\newcommand\R{\mathbb R}

\newcommand\tr{\operatorname{tr}}

\spnewtheorem*{main}{Main Theorem}{\bf}{\it}

\begin{document}
\title{Transformations on the set of all $n$-dimensional subspaces of a
Hilbert space preserving principal angles}
\titlerunning{Preserving principal angles between subspaces of a
Hilbert space}
\author{Lajos Moln\'ar}
\authorrunning{L. Moln\'ar}
\dedication{To my wife for her unlimited(?) patience}
\institute{Institute of Mathematics and Informatics,
           University of Debrecen,
           4010 Debrecen, P.O.Box 12, Hungary\\
\email{molnarl@math.klte.hu}}
\date{\today}
\communicated{H. Araki}
\maketitle
\begin{abstract}
Wigner's classical theorem on symmetry transformations plays
a fundamental role in quantum mechanics. It can be formulated, for
example,
in the following way: Every bijective transformation on the set $\L$ of
all 1-dimensional subspaces of a Hilbert space $H$ which
preserves
the angle between the elements of $\L$ is induced by either a unitary
or an antiunitary
operator on $H$. The aim of this paper is to extend Wigner's result from
the 1-dimensional case to the case of $n$-dimensional subspaces of
$H$ with $n\in \N$ fixed.
\end{abstract}

\section{Introduction and statement of the main result}

Let $H$ be a (real or complex) Hilbert space and denote $B(H)$ the
algebra of all bounded linear operators on $H$. By a projection we mean
a self-adjoint idempotent in $B(H)$.
For any $n\in \N$, $P_n(H)$ denotes the set of all
rank-$n$ projections on $H$, and $P_\infty(H)$ stands for the set of all
infinite rank projections.
Clearly, $P_n(H)$ can be identified with
the set of all $n$-dimensional subspaces of $H$. As it was mentioned in
the abstract, Wigner's theorem describes the bijective transformations
on the set $\L$ of all 1-dimensional subspaces of $H$ which preserve the
angle between the elements of $\L$. It seems to be a very natural
problem to try
to extend this result from the 1-dimensional case to the case of higher
dimensional subspaces (in our recent papers \cite{MolJMP},
\cite{MolCMP}, \cite{MolIJTP} we have presented several other
generalizations of Wigner's theorem for different structures). But what
about the angle between two
higher dimensional subspaces of $H$? For our present purposes, the most
adequate concept of
angles is that of the so-called principal angles (or canonical angles,
in a different terminology). This concept is
a generalization of the usual notion of angles between 1-dimensional
subspaces
and reads as follows. If $P,Q$ are finite dimensional projections, then
the principal angles between them (or, equivalently, between their
ranges as subspaces) is defined as the $\arccos$ of the square root of
the eigenvalues (counted according multiplicity) of the positive
(self-adjoint) finite
rank operator $QPQ$ (see, for example, \cite[Exercise VII.1.10]{Bha} or
\cite[Problem 559]{Kir}). We remark that this concept of angles
was motivated by the classical work \cite{Jor} of Jordan
and it has serious applications in statistics, for example (see the
canonical correlation theory of Hotelling \cite{Hot}, and also see the
introduction of \cite{Miao}).
The system of all principal angles between $P$ and $Q$ is denoted by
$\angle (P,Q)$.
Thus, we have the desired concept of angles between finite
rank projections. But in what follows we would like to extend Wigner's
theorem also for the case of infinite rank projections.
Therefore, we need the concept of principal angles also between infinite
rank projections.
Using deep concepts of
operator theory (like scalar-valued spectral measure and multiplicity
function) this could be carried out, but in order to formulate a
Wigner-type
result we need only the equality of angles. Hence, we can
avoid these complications saying that for arbitrary projections
$P,Q,P',Q'$ on $H$ we have $\angle (P,Q)=\angle
(P',Q')$ if and only if the positive operators $QPQ$ and $Q'P'Q'$ are
unitarily
equivalent. This obviously generalizes the equality of principal
angles between pairs of finite rank projection.

Keeping in mind the formulation of Wigner's theorem given in the
abstract,
we are now in a position to formulate the main result of the paper
which, we believe, also has physical interpretation.

\begin{main}
Let $n\in \N$. Let $H$ be a real or complex
Hilbert space with $\dim H\geq n$.
Suppose that $\phi :P_n(H)\to P_n(H)$ is a transformation
with the property that
\[
\angle (\phi(P),\phi(Q))=
\angle (P,Q)\qquad (P,Q\in P_n(H)).
\]
If $n=1$ or $n\neq \dim H/2$,
then there exists a linear or conjugate-linear isometry $V$ on $H$ such
that
\[
\phi(P)=VPV^* \qquad (P\in P_n(H)).
\]
If $H$ is infinite dimensional, the transformation $\phi:P_\infty (H)
\to P_\infty(H)$ satisfies
\[
\angle (\phi(P),\phi(Q))=
\angle (P,Q)\qquad (P,Q\in P_\infty (H)),
\]
and $\phi$ is surjective, then
there exists a unitary or antiunitary operator $U$ on $H$ such
that
\[
\phi(P)=UPU^* \qquad (P\in P_\infty (H)).
\]
\end{main}

As one can suspect from the formulation of our main result, there is a
system of exceptional cases, namely, when we have $\dim H=2n, n>1$.
In the next section we show that in those cases there do exist
transformations on $P_n(H)$ which preserve the principal angles but
cannot be written in the form appearing in our main theorem above.

\section{Proof}

This section is devoted to the proof of our main theorem.
In fact, this will follow from the statements below.

The idea of the proof can be summarized in a single sentence as follows.
We extend our transformation from $P_n(H)$ to a Jordan homomorphism
of the algebra $F(H)$ of all finite rank operators on $H$ which
preserves the rank-1 operators.
Fortunately, those maps turn to have a form and using this we can
achieve
the desired conclusion. On the other hand, quite unfortunately, we have
to work hard to carry out all the details of the proof that we
are just going to begin.

From now on, let $H$ be a real or complex Hilbert space and
let $n\in \N$. Since our statement obviously holds when $\dim H=n$,
hence we suppose that $\dim H>n$.

In the sequel, let $\tr$ be the usual trace functional on operators.
The ideal of all finite rank operators in $B(H)$ is denoted by $F(H)$.
Clearly, every element of $F(H)$ has a finite trace.
We denote by $F_s(H)$ the set of all self-adjoint elements of $F(H)$.

We begin with two key lemmas.
In order to understand why we consider the property \eqref{E:proji13}
in Lemma~\ref{L:proji1}, we note that if $\angle (P,Q)=\angle (P',Q')$
for some
finite rank projections $P,Q,P',Q'$, then, by definition, the positive
operators $QPQ$ and $Q'P'Q'$ are unitarily equivalent. This implies that
$\tr QPQ=\tr Q'P'Q'$. But, by the properties of the trace, we have $\tr
QPQ=\tr PQQ=\tr
PQ$ and, similarly, $\tr Q'P'Q'=\tr P'Q'$. So, if our transformation
preserves the principal angles between projections, then it necessarily
preserves the trace of the product of the projections in question.
This justifies the condition \eqref{E:proji13} in the next lemma.

\begin{lemma}\label{L:proji1}
Let $\P$ be any set of finite rank projections on $H$.
If $\phi :\P \to \P$ is a transformation with the property that
\begin{equation}\label{E:proji13}
\tr \phi(P)\phi(Q) =\tr PQ \qquad (P,Q\in \P),
\end{equation}
then $\phi$ has a unique real-linear extension $\Phi$ onto the
real-linear span $\bur_\R \P$ of $\P$. The transformation $\Phi$ is
injective, preserves the trace and satisfies
\begin{equation}\label{E:proji14}
\tr \Phi(A)\Phi(B) =\tr AB \qquad (A,B\in \bur_\R \P).
\end{equation}
\end{lemma}

\begin{proof}
For any finite sets $\{ \lambda_i \} \subset \R$ and $\{ P_i\} \subset
\P$ we define
\[
\Phi(\sum_i \lambda_i P_i)=
\sum_i \lambda_i \phi(P_i).
\]
We have to show that $\Phi$ is well-defined.
If $\sum_i \lambda_i P_i=\sum_k \mu_k Q_k$, where
$\{ \mu_k\} \subset \R$ and $\{ Q_k\} \subset \P$ are finite subsets,
then for any $R\in \P$ we compute
\[
\tr (\sum_i \lambda_i \phi(P_i)\phi(R))=
\sum_i \lambda_i \tr (\phi(P_i)\phi(R))=
\sum_i \lambda_i \tr (P_i R)=
\]
\[
\tr (\sum_i \lambda_i  P_i R)=
\tr (\sum_k \mu_k  Q_k R)=
\sum_k \mu_k  \tr(Q_k R)=
\]
\[
\sum_k \mu_k  \tr(\phi(Q_k)\phi(R))=
\tr(\sum_k \mu_k \phi(Q_k)\phi(R)).
\]
Therefore, we have
\[
\tr ((\sum_i \lambda_i \phi(P_i)-\sum_k \mu_k \phi(Q_k))\phi(R))=0
\]
for every $R\in \P$. By the linearity of the trace functional it follows
that
we have similar equality if we replace $\phi(R)$ by any finite linear
combination of $\phi(R)$'s. This gives us that
\[
\tr ((\sum_i \lambda_i \phi(P_i)-\sum_k \mu_k \phi(Q_k))
(\sum_i \lambda_i \phi(P_i)-\sum_k \mu_k \phi(Q_k)))=0.
\]
The operator
$(\sum_i \lambda_i \phi(P_i)-\sum_k \mu_k \phi(Q_k))^2$, being the
square of a self-adjoint operator, is positive. Since its trace is zero,
we obtain that
\[
(\sum_i \lambda_i \phi(P_i)-\sum_k \mu_k \phi(Q_k))^2=0
\]
which plainly implies that
\[
\sum_i \lambda_i \phi(P_i)-\sum_k \mu_k \phi(Q_k)=0.
\]
This shows that $\Phi$ is well-defined. The real-linearity of $\Phi$ now
follows from the definition.
The uniqueness of $\Phi$ is also trivial to see.
From \eqref{E:proji13} we immediately obtain \eqref{E:proji14}.
One can introduce an inner product on $F_s(H)$ by the formula
\[
\la A,B\ra =\tr AB \qquad (A,B\in F_s(H))
\]
(the norm induced by this inner product is called the Hilbert-Schmidt
norm).
The equality \eqref{E:proji14} shows that $\Phi$ is an isometry with
respect to this norm. Thus, $\Phi$ is injective. It follows from
\eqref{E:proji13} that
\[
\tr \phi(P)=\tr \phi(P)^2=\tr P^2=\tr P \qquad (P\in \P)
\]
which clearly implies that
\[
\tr \Phi(A)=\tr A \qquad (A\in \bur_\R \P).
\]
This completes the proof of the
lemma.
\qed
\end{proof}

In what follows we need the concept of Jordan homomorphisms. If $\cal
A$ and $\cal B$ are algebras, then a linear transformation $\Psi: \cal A
\to \cal B$ is called a Jordan homomorphism if it satisfies
\[
\Psi(A^2)=\Psi(A)^2 \qquad (A\in \cal A),
\]
or, equivalently, if
\[
\Psi(AB+BA)=\Psi(A)\Psi(B)+\Psi(B)\Psi(A) \qquad (A,B \in \cal A).
\]
Two projections
$P,Q$ on $H$ are said to be orthogonal if $PQ=QP=0$
(this means that the ranges of $P$ and $Q$ are orthogonal to each
other). In this case we write $P\perp Q$. We
denote $P\leq Q$ if $PQ=QP=P$
(this means that the range of $P$ is included in the range of $Q$).
In what follows, we shall use the following useful notation. If $x,y\in
H$, then $x\ot y$ stands for the operator defined by
\[
(x\ot y)z=\la z,y\ra x \qquad (z\in H).
\]

\begin{lemma}\label{L:proji2}
Let $\Phi:F_s(H) \to F_s(H)$ be a real-linear transformation which
preserves the rank-1 projections and the orthogonality between them.
Then there is an either linear or conjugate-linear
isometry $V$ on $H$ such that
\[
\Phi(A)=VAV^* \qquad (A\in F_s(H)).
\]
\end{lemma}

\begin{proof}
Since every finite-rank projection is the finite sum of pairwise
orthogonal rank-1 projections,
it is obvious that $\Phi$ preserves the finite-rank projections. It
follows from \cite[Remark 2.2]{BreSem} and the spectral theorem
that $\Phi$ is a Jordan homomorphism
(we note that \cite[Remark 2.2]{BreSem} is about self-adjoint operators
on finite dimensional complex Hilbert spaces, but the same argument
applies for $F_s(H)$ even if it is infinite dimensional and/or real).

We next prove that $\Phi$ can be extended to a Jordan
homomorphism of $F(H)$. To see this, first suppose that $H$ is complex
and consider the transformation
$\tilde \Phi: F(H) \to F(H)$ defined by
\[
\tilde \Phi(A+iB)=\Phi(A)+i\Phi(B) \qquad (A,B \in F_s(H)).
\]
It is easy to see that $\Phi$ is a linear
transformation which satisfies $\tilde \Phi(T^2)=\tilde
\Phi(T)^2$ $(T\in F(H))$.
This shows that $\tilde \Phi$ is a Jordan homomorphism.

If $H$ is real, then the situation is not so simple, but we can
apply a deep algebraic result of Martindale as follows (cf. the
proof of \cite[Theorem 3]{MolJAMS}). Consider the unitalized
algebra $F(H)\oplus \mathbb R I$ (of course, we have to add the identity
only when $H$ is infinite dimensional). Defining $\Phi(I)=I$,
we can extend $\Phi$ to
the set of all symmetric elements of the enlarged algebra in an obvious
way. Now we are in a position to apply the results in \cite{Mar}
on the extendability of Jordan
homomorphisms defined on the set of symmetric elements of a ring with
involution. To be precise, in \cite{Mar} Jordan homomorphism means an
additive map $\Psi$ which, besides $\Psi(s^2)=\Psi(s)^2$, also
satisfies $\Psi(sts)=\Psi(s)\Psi(t)\Psi(s)$. But if the ring in question
is 2-torsion free (in particular, if it is an algebra), this second
equality follows from the first one (see, for example, the proof of
\cite[6.3.2 Lemma]{Pal}).
The statements \cite[Theorem 1]{Mar} in the case when $\dim H
\geq 3$ and \cite[Theorem 2]{Mar} if $\dim H =2$ imply that
$\Phi$ can be uniquely extended to an associative homomorphism of
$F(H)\oplus \mathbb R I$ into itself. To be honest, since the
results of Martindale
concern rings and hence linearity does not appear,
we could guarantee only the additivity of the extension of $\Phi$.
However, the construction in \cite{Mar} shows that in the
case of algebras, linear Jordan homomorphisms have linear
extensions.

To sum up, in every case we have a Jordan homomorphism of $F(H)$
extending $\Phi$. In order to simplify the notation, we use the same
symbol $\Phi$ for the extension as well.

As $F(H)$ is a locally matrix ring
(every finite subset of $F(H)$ can be included in a subalgebra of $F(H)$
which is isomorphic to a full matrix algebra),
it follows from a classical result
of Jacobson and Rickart \cite[Theorem 8]{JR} that $ \Phi$ can
be written as $\Phi= \Phi_1 +\Phi_2$, where $\Phi_1$ is a
homomorphism and $\Phi_2$ is an antihomomorphism. Let $P$ be a
rank-1 projection on $H$. Since $\Phi(P)$ is also rank-1,
we obtain that
one of the idempotents $\Phi_1(P), \Phi_2(P)$ is zero. Since $F(H)$ is a
simple ring, it is easy to see that this implies that either
$\Phi_1$ or $\Phi_2$ is identically zero, that is, $ \Phi$ is
either a homomorphism or an antihomomorphism of $F(H)$. In what follows
we can assume without loss of generality that $\Phi$ is a
homomorphism. Since the kernel of $\Phi$ is an ideal in $F(H)$ and
$F(H)$ is simple, we obtain that $\Phi$ is injective.

We show that $\Phi$ preserves the rank-1 operators.
Let $A\in F(H)$ be of rank 1.
Then there is a rank-$1$ projection $P$ such that
$PA=A$.
We have $\Phi(A)= \Phi(PA)=\Phi(P) \Phi(A)$
which proves that $\Phi(A)$ is of rank at most $1$.
Since $\Phi$ is injective, we obtain that the rank of $\Phi(A)$ is
exactly 1.
From the conditions of the lemma it follows that $\phi$ sends rank-2
projections to rank-2 projections.
Therefore, the range of $\Phi$ contains an operator with rank greater
than 1. We now refer to Hou's work \cite{Hou} on the form of linear
rank preservers on operator algebras.
It follows from the argument leading to \cite[Theorem 1.3]{Hou}
that either there are linear
operators $T,S$ on $H$ such that $\Phi$ is of the form
\begin{equation*}
\Phi(x\ot y)=(Tx)\ot (Sy) \qquad (x,y \in H)
\end{equation*}
or there are conjugate-linear operators $T',S'$ on $H$ such that $\Phi$
is of the form
\begin{equation}\label{E:proji12b}
\Phi(x\ot y)=(S'y)\ot (T'x) \qquad (x,y \in H).
\end{equation}
Suppose that we have the first possibility. By the multiplicativity of
$\Phi$ we obtain that
\begin{equation}\label{E:proji18}
\begin{gathered}
\la u,y\ra Tx\ot Sv=
\la u,y\ra\Phi(x\ot v)=
\Phi(x\ot y \cdot u\ot v)=
\\
\Phi(x\ot y)\Phi(u\ot v)=
\la Tu, Sy\ra Tx\ot Sv.
\end{gathered}
\end{equation}
This gives us that $\la Tu, Sy\ra=\la u,y\ra$ for every $u,y\in H$. On
the other hand, since $\Phi$ sends rank-1 projections to
rank-1 projections,
we obtain that for every unit vector $x\in H$ we have $Tx=Sx$.
These imply that $T=S$ is an isometry and with the notation $V=T=S$ we
have
\[
\Phi(A)=VAV^*
\]
for every $A\in F_s(H)$.

We show that the possibility \eqref{E:proji12b} cannot occur. In fact,
similarly to \eqref{E:proji18} we have
\[
\la u,y\ra S'v\ot T'x=
\la S'v,T'x\ra S'y\ot T'u \qquad (x,y,u,v \in H).
\]
Fixing unit vectors $x=y=u$ in $H$ and considering the operators above
at $T'x$, we find that
\[
S'v=\la S'v,T'x\ra\la T'x,T'u\ra S'y
\]
giving us that $S'$ is of rank 1. Since $\Phi$ sends rank-2 projections
to rank-2 projections, we arrive at a contradiction. This completes the
proof of the lemma.
\qed
\end{proof}

We are now in a position to present a new proof of the nonsurjective
version of Wigner's theorem which is equivalent to the statement of our
main theorem in the case when $n=1$. For another proof see \cite{Sha}.

To begin, observe that if $P,Q$ are finite rank projections such
that $\tr PQ=0$,
then we have $\tr (PQ)^*PQ=\tr QPQ=\tr PQQ=\tr PQ=0$ which implies that
$(PQ)^*(PQ)=0$. This gives us that $PQ=0=QP$. Therefore, $P$ is
orthogonal to $Q$ if and only if $\tr PQ=0$.

\begin{theorem}\label{T:proji1}
Let $\phi:P_1(H) \to P_1(H)$ be a transformation with the property that
\begin{equation}\label{E:proji15}
\tr \phi(P)\phi(Q)=\tr PQ \qquad (P, Q\in P_1(H)).
\end{equation}
Then there is an either linear or conjugate-linear isometry $V$ on $H$
such that
\[
\phi(P)=VPV^* \qquad (P\in P_1(H)).
\]
\end{theorem}

\begin{proof}
By the spectral theorem it is obvious that the real linear span of
$P_1(H)$ is $F_s(H)$.
Then, by Lemma~\ref{L:proji1} we see that there is a unique real-linear
extension $\Phi$ of $\phi$ onto $F_s(H)$
which preserves the
rank-1 projections and, by \eqref{E:proji15},
$\Phi$ also preserves the orthogonality between the elements of
$P_1(H)$. Lemma~\ref{L:proji2} applies to complete the proof.
\qed
\end{proof}

As for the cases when $n>1$ we need the following lemma.
Recall that we have previously supposed that $\dim H>n$.

\begin{lemma}\label{L:proji3}
Let $1<n \in \N$. Then $\bur_\R P_n(H)$ coincides with $F_s(H)$.
\end{lemma}

\begin{proof}
Since the real-linear span of $P_1(H)$ is $F_s(H)$, it
is sufficient to show that every rank-1 projection is a real-linear
combination of rank-$n$ projections.
To see this, choose orthonormal vectors $e_1, \ldots, e_{n+1}$ in $H$.
Let $E=e_1\ot e_1 +\ldots +e_{n+1}\ot e_{n+1}$ and
define
\begin{equation*}
P_k=E-e_k\ot e_k \qquad (k=1,\ldots, n+1).
\end{equation*}
Clearly, every $P_k$
can be represented by a $(n+1)\times (n+1)$ diagonal matrix whose
diagonal entries are all 1's with the exception of the $k^{\text{th}}$
one which is 0.
The equation
\[
\lambda_1 P_1 + \ldots +\lambda_{n+1} P_{n+1}=e_1\ot e_1
\]
gives rise to a system of linear equations with unknown scalars
$\lambda_1 , \ldots, \lambda_{n+1}$.
The matrix of this system of equations is an
$(n+1)\times (n+1)$ matrix whose diagonal consists of $0$'s
and its off-diagonal entries are all 1's.
It is easy to see that this matrix is nonsingular, and hence
$e_1\ot e_1$ (and, similarly, every other $e_k\ot e_k$) is a real-linear
combination of $P_1, \ldots , P_{n+1}$. This completes the proof.
\qed
\end{proof}

We continue with a technical lemma.

\begin{lemma}\label{L:proji4}
Let $P,Q$ be projections on $H$. If $QPQ$ is a projection, then
there are pairwise orthogonal projections $R,R',R''$ such that
$P=R+R'$, $Q=R+R''$.
In particular, we obtain that $QPQ$ is a projection if and only
if $PQ=QP$.
\end{lemma}

\begin{proof}
Let $R=QPQ$.
Since $R$ is a projection whose range is contained in the range of $Q$,
it follows that $R''=Q-R$ is a
projection which is orthogonal to $R$.

If $x$ is a unit vector in the range of $R$, then we have
$\|QPQx\|=1$. Since $PQx$ is a vector whose norm is at most 1 and its
image under the projection $Q$ has norm 1, we obtain that $PQx$ is a
unit
vector in the range of $Q$. Similarly, we obtain that $Qx$ is a unit
vector in the range of $P$ and, finally, that $x$ is a unit vector
in the range of $Q$. Therefore, $x$ belongs to the range of $P$ and
$Q$.
Since $x$ was arbitrary, we can infer that the range of $R$ is included
in the range of $P$. Thus, we obtain
that $R'=P-R$ is a projection which is orthogonal to $R$.

Next, using the obvious relations
\[
PR=RP=R,\quad  QR=RQ=R
\]
we deduce
\begin{equation}\label{E:proji4}
\begin{gathered}
(Q-R)(P-R)(Q-R)=\\
QPQ-QPR-QRQ+QR-RPQ+RPR+RQ-R=\\
R-R-R+R-R+R+R-R=0.
\end{gathered}
\end{equation}
Since $A^*A=0$ implies $A=0$ for any $A\in B(H)$, we obtain from
\eqref{E:proji4} that $R'R''=(P-R)(Q-R)=0$.

The second part of the assertion is now easy to check.
\qed
\end{proof}

We next prove the assertion of our main theorem in the case when $1<n\in
\N$ and $H$ is infinite dimensional.

\begin{theorem}\label{T:proji2}
Suppose $1<n\in \N$ and $H$ is infinite dimensional.
If $\phi: P_n(H) \to P_n(H)$ is a transformation such that
\[
\angle (\phi(P),\phi(Q))=
\angle (P,Q)\qquad (P,Q\in P_n(H)),
\]
then there exists a linear or conjugate-linear isometry $V$ on $H$ such
that
\begin{equation*}
\phi(P)=VPV^* \qquad (P\in P_n(H)).
\end{equation*}
\end{theorem}

\begin{proof}
By Lemma~\ref{L:proji1} and Lemma~\ref{L:proji3}, $\phi$ can be uniquely
extended to an injective real-linear
transformation $\Phi$ on $F_s(H)$. The main point of the proof is to
show that $\Phi$ preserves the rank-1 projections.
In order to verify this,
just as in the proof of Lemma~\ref{L:proji3}, we consider
orthonormal vectors $e_1, \ldots, e_{n+1}$ in $H$, define
$E=e_1\ot e_1 +\ldots +e_{n+1}\ot e_{n+1}$ and set
\begin{equation*}
P_k=E-e_k\ot e_k \qquad (k=1,\ldots, n+1).
\end{equation*}

We show that the ranges of all $P_k'=\phi(P_k)$'s can be jointly
included in an $(n+1)$-dimensional subspace of $H$.
To see this, we first recall that $\Phi$ has the property that
\begin{equation*}
\tr \Phi(A)\Phi(B)=\tr AB \qquad (A,B\in F_s(H))
\end{equation*}
(see Lemma~\ref{L:proji1}).
Next we have the following property of $\Phi$: if $P,Q$ are
orthogonal rank-1 projections, then $\Phi(P)\Phi(Q)=0$. Indeed,
if $P,Q$  are orthogonal, then we can include them into two orthogonal
rank-$(n+1)$ projections. Now, referring to the construction given in
Lemma~\ref{L:proji3} and having in mind that $\Phi$ preserves the
orthogonality between rank-$n$ projections, we obtain that
$\Phi(P)\Phi(Q)=0$. (Clearly, the same argument works if $\dim H\geq
2(n+1)$.)
Since the rank-$n$ projections $P_k$ are commuting, by the preserving
property of $\phi$ and Lemma~\ref{L:proji4}, it follows that the
projections $\Phi(P_k)$ are also commuting.
It is well-known that any finite commuting family of operators in
$F_s(H)$ can be diagonalized by the same unitary transformation (or, in
the real case, by the same orthogonal transformation).
Therefore, if we resctrict $\Phi$ onto the real-linear subspace in
$F_s(H)$
generated by $P_1, \ldots, P_{n+1}$, then it can be identified with a
real-linear operator from $\R^{n+1}$ to $\R^m$ for some $m\in \N$.
Clearly, this restriction of $\Phi$ can be represented by
an $m\times (n+1)$ real matrix $T=(t_{ij})$. Let us examine how the
properties of $\Phi$ are reflected in those of the matrix $T$.
First, $\Phi$ is trace preserving. This gives us that for every
${\underline \lambda}\in \R^{n+1}$
the sums of the coordinates of the vectors $T {\underline\lambda}$
and $\underline \lambda$ are the same.
This easily implies that the sum of the entries of $T$ lying in a fixed
column is always 1.
As we have already noted,
$\Phi(e_i\ot e_i)\Phi(e_j\ot e_j)=0$ holds
for every $i\neq j$. For the matrix $T$ this means
that the coordinatewise product of any two
columns of $T$ is zero. Consequently in every row of $T$ there is at
most one nonzero entry.
Since $\Phi$ sends rank-$n$ projections to rank-$n$ projections, we see
that this possibly nonzero entry is necessarily 1.
So, every row contains
at most one 1 and all the other entries in that row are 0's.
Since the sum of the elements in every column is 1, we have that in
every column there is exactly one 1 and all
the other entries are 0's in that column.
These now easily imply that if
${\underline\lambda}\in \R^{n+1}$ is such that its coordinates
are all 0's with the exception of one which is 1,
then $T{\underline \lambda}$ is of the same kind.
What concerns $\Phi$, this means that $\Phi$ sends every $e_k\ot e_k$
$(k=1, \ldots, n+1)$ to a rank-1 projection.

So, we obtain that
$\Phi$ preserves the rank-1 projections and the orthogonality between
them.
Now, by Lemma~\ref{L:proji2} we conclude the proof.
\qed
\end{proof}

We turn to the case when $H$ is finite dimensional.

\begin{theorem}\label{T:proji5}
Suppose $1<n \in \N$, $H$ is finite dimensional and $n\neq
\dim H/2$. If $\phi :P_n(H) \to P_n(H)$ satisfies
\[
\angle (\phi(P),\phi(Q))=
\angle (P,Q)\qquad (P,Q\in P_n (H)),
\]
then there exists a unitary or antiunitary operator $U$ on $H$ such
that
\begin{equation}\label{E:proji30}
\phi(P)=UPU^* \qquad (P\in P_\infty (H)).
\end{equation}
\end{theorem}

\begin{proof}
First suppose that $\dim H=2d$, $1<d\in \N$. If $n=1, \ldots, d-1$, then
we can apply the method followed in the proof of Theorem~\ref{T:proji2}
concerning the infinite dimensional case.
If $n=d+1, \ldots, 2d-1$, then consider the transformation
$\psi: P\mapsto I-\phi(I-P)$ on $P_{2d-n}(H)$.
We learn from \cite[Problem 559]{Kir} that if $\angle
(P,Q)=\angle (P',Q')$, then there exists a unitary operator $U$ such
that $UPU^*=P'$ and $UQU^*=Q'$. It follows from the preserving
property of $\phi$ that for any $P,Q \in P_{2d-n}(H)$ we have
\[
\phi(I-P)=U(I-P)U^*, \quad
\phi(I-Q)=U(I-Q)U^*
\]
for some unitary operator $U$ on $H$.
This gives us that
\[
\angle( \psi(P), \psi(Q))=
\angle( UPU^*, UQU^*)=
\angle( P, Q).
\]
In that way we can reduce the problem to the previous
case. So, there is an either unitary or antiunitary operator $U$ on $H$
such that
\[
\psi(P)=UPU^* \qquad (P\in P_{2d-n}(H)).
\]
It follows that $\phi(I-P)=I-\psi(P)=I-UPU^*=U(I-P)U^*$, and hence we
have the result for the considered case.

Next suppose that $\dim H=2d+1$, $d\in \N$. If $n=1, \ldots, d-1$,
then once again we
can apply the method followed in the proof of Theorem~\ref{T:proji2}.
If $n=d+2, \ldots, 2d+1$, then using the 'dual method'
that we have applied right above we can reduce the problem to the
previous case.
If $n=d$, consider a fixed rank-$d$ projection $P_0$. Clearly, if $P$ is
any rank-$d$ projection orthogonal to $P_0$, then the rank-$d$
projection $\phi(P)$ is orthogonal to
$\phi(P_0)$. Therefore, $\phi$ induces a transformation $\phi_0$ between
$d+1$-dimensional spaces (namely, between the orthogonal complement of
the
range of $P_0$ and that of the range of $\phi(P_0)$) which preserves the
principal angles between
the rank-$d$ projections. Our 'dual method' and the result concerning
1-dimensional subspaces lead us to the conslusion
that the linear extension of $\phi_0$ maps rank-1 projections to rank-1
projections and preserves the orthogonality between them. This implies
that the same holds true for our original
transformation $\phi$. Just as before, using Lemma~\ref{L:proji1} and
Lemma~\ref{L:proji2} we can conclude the proof.
In the remaining case $n=d+1$ we apply the 'dual method' once again.
\qed
\end{proof}

We now show that the case when $1<n \in \N$, $n=\dim H/2$ is really
exceptional. To see this, consider the transformation
$\phi: P \mapsto I-P$ on $P_n(H)$. This maps $P_n(H)$ into itself and
preserves the principal angles. As for the complex case, the preserving
property follows
from \cite[Exercise VII.1.11]{Bha} while in the real case it was proved
already by Jordan in \cite{Jor} (see \cite{Paige}, p. 310).
Let us now suppose that the transformation $\phi$ can be written in the
form \eqref{E:proji30}.
Pick a rank-1 projection $Q$ on $H$. We know that it is a real linear
combination of some $P_1, \ldots, P_{n+1}\in P_n(H)$. It would follow
from \eqref{E:proji30}
that considering the same linear combination of
$\phi(P_1), \ldots, \phi(P_{n+1})$, it is a rank-1 projection as well.
But due to the definition of $\phi$, we get that this linear combination
is a constant minus $Q$. By the trace preserving property we
obtain that this constant is $1/n$. Since $n>1$, the operator
$(1/n)I-Q$ is obviously not a projection.
Therefore, we have arrived at a contradiction. This shows that the
transformation above can not be written in the form \eqref{E:proji30}.

It would be a nice result if one could prove that in the present case
(i.e., when $1<n, n=\dim H/2$) up to unitary-antiunitary equivalence,
there are exactly two transformations on $P_n(H)$ preserving principal
angles, namely, $P\mapsto P$ and $P\mapsto I-P$. This is left as an open
problem.

We now turn to our statement concerning infinite rank
projections. In the proof we shall use the following simple lemma.
If $A\in B(H)$, then denote by $\rng A$ the range of $A$.

\begin{lemma}\label{L:proji5}
Let $H$ be an infinite dimensional Hilbert space.
Suppose $P,Q$ are projections on $H$ with the property that for any
projection $R$ with finite corank we have $RP=PR$ if and only if
$RQ=QR$. Then either $P=Q$ or $P=I-Q$.
\end{lemma}

\begin{proof}
Let $R$ be any projection on $H$ commuting with $P$. By
Lemma~\ref{L:proji4}, it is easy to see that we can
choose a monotone decreasing net $(R_\alpha)$ of projections with finite
corank such that $(R_\alpha)$ converges weakly to $R$ and
$R_\alpha$ commutes
with $P$ for every $\alpha$. Since $R_\alpha$ commutes with $Q$ for
every $\alpha$,
we obtain that $R$ commutes with $Q$. Interchanging the role of $P$ and
$Q$, we obtain that any projection commutes with $P$ if and only if it
commutes with $Q$.

Let $x$ be any unit vector from the range of $P$. Consider $R=x\ot x$.
Since $R$ commutes with $P$, it must commute with $Q$ as well. By
Lemma~\ref{L:proji4} we obtain that $x$ belongs either to the range of
$Q$ or to its orthogonal complement. It follows that
either $d(x,\rng Q)=0$, or $d(x, \rng Q)=1$. Since the set of all unit
vectors in the range of $P$ is connected and the distance function is
continuous, we get that either every unit vector in $\rng P$ belongs to
$\rng Q$ or every unit vector in $\rng P$ belongs to $(\rng Q)^\perp$.
Interchanging the role of $P$ and $Q$, we find that
either $\rng P=\rng Q$ or $\rng P=(\rng Q)^\perp$.
This gives us that either $P=Q$ or $P=I-Q$.
\qed
\end{proof}

\begin{theorem}
Let $H$ be an infinite dimensional Hilbert space. Suppose that
$\phi: P_\infty(H)\to P_\infty(H)$ is a surjective transformation with
the property that
\[
\angle (\phi(P),\phi(Q))=
\angle (P,Q)\qquad (P,Q\in P_\infty (H)).
\]
Then there exists a unitary or antiunitary operator $U$ on $H$ such that
\[
\phi(P)=UPU^* \qquad (P\in P_\infty (H)).
\]
\end{theorem}

\begin{proof}
We first prove that $\phi$ is injective.
If $P,P' \in P_\infty (H)$ and $\phi(P)=\phi(P')$, then by the
preserving property of $\phi$ we have
\begin{equation}\label{E:proji10}
\angle (P,Q)=\angle (P', Q) \qquad (Q\in P_\infty (H)).
\end{equation}
Putting $Q=I$, we see that $P$ is unitarily equivalent to $P'$.
We distiguish two cases.
First, let $P$ be of infinite corank.
By \eqref{E:proji10}, we deduce that for every $Q\in P_\infty
(H)$ we have $Q\perp P$ if and only if $Q\perp P'$. This gives us that
$P=P'$.
As the second possibility, let $P$ be of finite corank. Then $P,P'$ can
be written in
the form $P=I-P_0$ and $P'=I-P_0'$, where, by the equivalence of $P,P'$,
the projections $P_0$ and $P_0'$ have finite and equal rank.
Let $Q_0$ be any finite rank projection on $H$. It follows from
\[
\angle (I-P_0, I-Q_0)=
\angle (I-P_0', I-Q_0)
\]
that there is a unitary operator $W$ on $H$ such that
\[
W(I-Q_0)(I-P_0)(I-Q_0)W^*=
(I-Q_0)(I-P_0')(I-Q_0).
\]
This implies that
\begin{equation*}
\begin{gathered}
W(-Q_0-P_0+P_0Q_0+Q_0P_0-Q_0P_0Q_0)W^*=\\
-Q_0-P_0'+P_0'Q_0+Q_0P_0'-Q_0P_0'Q_0.
\end{gathered}
\end{equation*}
Taking traces, by the equality of the rank of $P_0$ and $P_0'$, we
obtain that
\begin{equation}\label{E:proji19}
\tr P_0Q_0=\tr P_0'Q_0.
\end{equation}
Since this holds for every finite rank projection $Q_0$ on $H$, it
follows
that $P_0=P_0'$ and hence we have $P=P'$. This proves the injectivity of
$\phi$.

Let $P \in P_\infty (H)$ be of infinite corank. Then there is a
projection $Q\in P_\infty (H)$ such that $Q\perp P$. By the preserving
property of $\phi$, this implies that $\phi(Q)\perp \phi(P)$ which means
that $\phi(P)$ is of infinite corank. One can similarly prove that if
$\phi(P)$ is of infinite corank, then the same must hold for $P$.
This yields that $P\in P_\infty (H)$ is of finite corank if and only if
so is $\phi(P)$.

Denote by $P_f(H)$ the set of all finite rank projections on $H$.
It follows that the transformation $\psi: P_f(H)\to P_f(H)$
defined by
\[
\psi(P)=I-\phi(I-P) \qquad (P\in P_f(H))
\]
is well-defined and bijective.
Since $\phi(I-P)$ is unitarily equivalent to $I-P$ for every $P\in
P_f(H)$
(this is because
$\angle (\phi(I-P), \phi(I-P))= \angle (I-P, I-P)$),
it follows that $\psi$
is rank preserving.

We next show that
\begin{equation}\label{E:proji11}
\tr \psi(P)\psi(Q)=\tr PQ \qquad (P,Q\in P_f(H)).
\end{equation}
This can be done following the argument leading to \eqref{E:proji19}.
In fact, by the preserving property of $\phi$ there is a unitary
operator $W$ on $H$ such that
\[
W(I-\psi(Q))(I-\psi(P))(I-\psi(Q))W^*=
(I-Q)(I-P)(I-Q).
\]
This gives us that
\begin{equation*}
\begin{gathered}
W(-\psi(Q)-\psi(P)+\psi(P)\psi(Q)+\psi(Q)\psi(P)-\psi(Q)\psi(P)\psi(Q))W^*=\\
-Q-P+PQ+QP-QPQ.
\end{gathered}
\end{equation*}
Taking traces on both sides and referring to the rank preserving
property of $\psi$, we obtain \eqref{E:proji11}.
According to Lemma~\ref{L:proji1}, let $\Psi:F_s(H)\to F_s(H)$ denote
the unique real-linear extension of $\psi$ onto $\bur_\R
P_f(H)=F_s(H)$. We know that $\Psi$ is injective.
Since $P_f(H)$ is in the range of $\Psi$, we obtain that $\Psi$
is surjective as well. It is easy to see that Lemma~\ref{L:proji2}
can be applied and we
infer that there exists an either unitary or
antiunitary operator $U$ on $H$ such that
\[
\Psi(A)=UAU^* \qquad (A\in F_s(H)).
\]

Therefore, we have
\[
\phi(P)=UPU^*
\]
for every projection $P\in P_\infty(H)$ with finite corank. It remains
to prove that
the same holds true for every $P\in P_\infty(H)$ with infinite corank as
well. This could be
quite easy to show if we know that $\phi$ preserves the order between
the elements of $P_\infty(H)$. But this
property is far away from being easy to verify. So we choose a different
approach to attack the problem.

Let $P\in P_\infty(H)$ be a projection of infinite corank.
By the preserving property of $\phi$ we
see that for every $Q\in P_\infty (H)$ the operator
$\phi(Q)\phi(P)\phi(Q)$ is a projection if and if $QPQ$ is a projection.
By Lemma~\ref{L:proji4}, this means that
$\phi(Q)$ commutes with $\phi(P)$ if and only if $Q$ commutes with $P$.
Therefore, for any $Q\in P_\infty (H)$ of finite corank, we obtain that
$Q$ commutes with $U^*\phi(P)U$ (this is equivalent to that
$\phi(Q)=UQU^*$ commutes with $\phi(P)$) if and only if $Q$ commutes
with $P$.

By Lemma~\ref{L:proji5} we have two possibilities, namely,
either $U^*\phi(P)U=P$ or $U^*\phi(P)U=I-P$.
Suppose that
$U^*\phi(P)U=I-P$.
Consider a complete orthonormal basis
$e_0, e_\gamma \, ( \gamma \in \Gamma)$ in the range of $P$ and,
similarly, choose a complete orthonormal basis
$f_0, f_\delta \, (\delta\in \Delta)$ in the range of
$I-P$. Pick nonzero scalars $\lambda$, $\mu$ with the property
that $|\lambda |^2 +|\mu|^2=1$ and $|\lambda| \neq |\mu|$. Define
\[
Q=(\lambda e_0+\mu f_0)\ot (\lambda e_0+\mu f_0) +
\sum_{\gamma} e_\gamma \ot e_\gamma +
\sum_{\delta} f_\delta \ot f_\delta.
\]
Clearly, $Q$ is of finite corank (in fact, its corank is 1). Since
$\phi(Q)\phi(P)\phi(Q)=UQU^*\phi(P)UQU^*$ is unitarily
equivalent to $QPQ$, it
follows that the spectrum of
$QU^*\phi(P)UQ$ is equal to the
spectrum of $QPQ$. This gives us that the spectrum of
$Q(I-P)Q$ is equal to the spectrum of $QPQ$. By the construction of
$Q$ this means that
\[
\{ 0, 1, |\mu| ^2\} =
\{ 0, 1, |\lambda | ^2\}
\]
which is an obvious contradiction. Consequently, we have
$U^*\phi(P)U=P$, that is, $\phi(P)=UPU^*$.
Thus, we have proved that this latter equality holds for every $P\in
P_\infty (H)$ and the proof is complete.
\qed
\end{proof}

\begin{acknowledgements}
          This research was supported from the following sources:
          (1) Hungarian National Foundation for Scientific Research
          (OTKA), Grant No. T030082, T031995,
          (2) A grant from the Ministry of Education, Hungary, Reg.
          No. FKFP 0349/2000.
\end{acknowledgements}


\end{document}